\documentclass[11pt,a4paper]{article}

\usepackage{amssymb,latexsym}

\newtheorem{theorem}{Theorem}

\newtheorem{proposition}[theorem]{Proposition}

\newtheorem{lemma}[theorem]{Lemma}

\begin{document}

\title{On the Height Profile of a Conditioned \newline Galton-Watson Tree 
\thanks{partially supported by the German Research Foundation DFG}}
\author{G\"otz Kersting\\ University of Frankfurt am Main\thanks{Fachbereich 
Mathematik, Postfach 11 19 32, D-60054 Frankfurt/Main}}
\date{March 1998}
\maketitle

\begin{abstract}
Recently Drmota and Gittenberger (1997) proved a conjecture due to 
Aldous (1991) on the height profile of a Galton-Watson tree with an 
offspring distribution of finite variance, conditioned on a total size of $n$ 
individuals. The conjecture states that in distribution its shape, 
more precisely its scaled height profile 
coincides asymptotically with the local time process of a Brownian excursion 
of duration 1. We give a proof of the result, which extends to the 
case of an infinite variance offspring distribution. This requires a 
different strategy, since in the infinite variance case there is no 
longer a relationship to the local time of Brownian resp. L\'{e}vy excursions.
\end{abstract}

\small
AMS classification numbers: Primary 60J80, Secondary 60F17 
   
Key words: Galton-Watson trees, L\'{e}vy-excursions, functional limit 
theorems
\normalsize

\section{Introduction and main result} 

In this paper we analyse the shape of a Galton-Watson tree, 
conditioned to have a total number of $n$ individuals. More precisely 
we study its asymptotic height profile, as $n \rightarrow 
\infty$. 

By a tree $t$ we mean a rooted, ordered,  finite tree.  We consider 
the vertices to represent individuals, such that $t$ can be regarded 
as the family tree of the progeny of some founding ancestor (the 
root).  It is assumed that among siblings there is an order (of 
birth), which allows to imbed such trees into the plane.  The total 
number of vertices of $t$, its \emph{size} is denoted by $s(t)$.  We 
shall be interested in the way, in which the size of the
different generations vary along the tree.  An individual $i$ belongs to 
the $k$'th generation, if the path from $i$ to the root contains 
exactly $k$ edges.  Let $z_{k}$ denote the number of individuals in 
generation $k=0,1, \ldots$.  The sequence $ 1 = z_{0}, z_{1}, z_{2}, \ldots $ 
is called the  \emph{height profile} of $t$ (possibly after a 
suitable renormalisation).

For a Galton-Watson tree $T$ the number of children of 
the different individuals are assumed to be independent and 
identically distributed random variables. $p_{x},  x= 
0,1, \ldots$ denotes the probability, that an individual possesses $x$ 
children, and $Z_{k}$ is the number of individuals in 
generation $k$. We shall consider $T$, conditioned on a size
$s(T) = n$. Such a conditioned Galton-Watson tree will be abbreviated 
as $CGW(n)$-tree. Since it contains n individuals, 
it is natural to consider its height profile
\begin{displaymath}
	H_{u}^{n} = a_{n}^{-1} Z_{[nu/a_{n}]},\; u \geq 0,
\end{displaymath}
rescaled with some sequence $(a_{n})$ of positive numbers. The aim is to choose 
$(a_{n})$ such, that $H^{n}$ converges in distribution. Then $a_{n}$ 
gives the magnitude of the breadth of the $CGW(n)$-tree, whereas 
$n/a_{n}$ indicates the order of $h(T) = \max \{k \: | \: Z_{k} > 
0\}$, the height of the $CGW(n)$-tree. -  We assume:\bigskip 

\textbf{Assumption A} {\itshape \\ 
1.) The offspring distribution $(p_{x})_{x}$ has mean 1:
\begin{displaymath}
	\sum_{x=0}^{\infty} x p_{x} = 1.
\end{displaymath} 
The greatest common divisor of all $x$ with $p_{x} >0$ is 1.\\ 
2.) There are positive numbers $a_{n}$ such that $a_{n}^{-1}
(\xi_{1}+\ldots+\xi_{n}-n)$ converges in distribution to a 
non-degenerate limit law $\nu$, as $n \rightarrow 
\infty$. Here $\xi_{1}, \xi_{2}, \ldots$ 
denote independent random variables with distribution $(p_{x})_{x}$.}
\bigskip 

The case of an offspring distribution with a finite mean can be treated 
in much the same way and is not a genuine generalisation (see Kennedy 
\cite{ken}). Similar 
the assumption on the g.c.d. may be removed. Criticality and a g.c.d. 
equal to 1 are assumed just for convenience. 

The second assumption has been completely analysed. 
It has several implications, which are discussed in 
chapter XVII in Feller \cite{fe2} and in the monograph of Gnedenko 
and Kolmogorov \cite{gn}. In particular $(p_{x})$ is in 
the domain of attraction of a stable law. Its index $\alpha$ belongs 
to $(1,2]$, since we deal with an offspring distribution of
finite mean. This implies $a_{n} =o(n)$. The limit law $\nu = 
\nu_{\alpha}$ is determined by $\alpha$ up to a scaling constant. 
For $\alpha = 2$ this is simply the normal law, and 
for $1 < \alpha < 2$ it is a one-sided stable law, 
because the negative tail of the offspring distribution is zero.

We shall prove, that under these assumptions the height profile 
$H^{n}$ converges in distribution, where the $a_{n}$ are just the 
numbers, given in the assumption. The limiting process 
$H=(H_{u})_{u}$ turns out to be a functional of certain normalized 
excursions:  $H=\psi (Y)$. By a \emph{normalized excursion} we understand
a process $Y = (Y_{s})_{0 \leq s \leq 1}$, such 
that $Y_{0}=Y_{1}=0$ and $\inf_{\delta  \leq s \leq 1- 
\delta }  Y_{s} > 0$ for all $\delta >0$. 
$H$ will be obtained from $Y$ as follows: The corresponding 
\emph{cumulative height profile}
\begin{displaymath}
	C_{u} = \int_{0}^{u} H_{v} \: dv, \; u \geq 0,
\end{displaymath}
is given by
\begin{displaymath}
	C_{u}= \sup \{s \leq 1 \;: \; \int_{0}^{s} \frac{dt}{Y_{t}}  \leq u \}.
\end{displaymath}
$H$, having a.s. paths continuous from the right, is completely 
determined by the cumulative height process $C$. If $\int_{0+} dt/Y_{t}  
< \infty$, then $h>0$ and $H_{u} >0$ for all $u \in (0,h)$, also $C_{u} 
\rightarrow 1$, as $u \rightarrow h$. Then we say, that $H$ is 
\emph{non-degenerate}. The case 
$\int_{0+} dt/Y_{t}  = \infty$, implying $C \equiv 0$ and $H \equiv 
0$, is in the sequel of no significance. - The same transformation 
already appears in the paper  \cite{la} by Lamperti, who used it to 
transform L\'{e}vy-processes into continuous branching processes. Our 
theorem is in a way a conditioned version of Lamperti's result.

The normalized excursions are obtained as follows. There is
a unique L\'{e}vy-process $X=(X_{s})_{s \geq 0}$ (a process with 
independent and stationary increments), such that $X_{0}=0$ and the law 
$\nu$ is the distribution of $X_{1}$, and to each of these processes
there belongs a normalized excursion $Y$.  This is explained in detail 
in Bertoin's monograph \cite{ber} (in particular chapter VIII.4).

\begin{theorem} Let $H^{n}$ be the scaled height profile of a 
$CGW(n)$-tree, satisfying assumption A. Then, as $n \rightarrow 
\infty$, $H^{n}$ converges in distribution to the process 
$H = \psi (Y)$, derived from the corresponding normalized 
L\'{e}vy-excursion $Y$. $H$ is a.s. non-degenerate.
\end{theorem}

In general the excursions contain jumps. These jumps also show 
up in $H$. Thus we regard $H^{n}$ and $H$ as random elements in the  
space of c\'{a}dl\'{a}g-functions, endowed with the usual Skorohod  
$J_{1}$ topology (compare  \cite{et}). 

A neat case is that of an offspring distribution with finite 
variance:
\begin{displaymath}
	\sigma^{2} = \sum_{x=0}^{\infty} x^{2} p_{x} \; - 1 < \infty .
\end{displaymath}
Then, choosing $a_{n} = \sigma  n^{1/2}$, $Y$ is a 
normalized Brownian excursion. In this situation the limiting process 
allows another appealing description: $(H_{u})_{u} \stackrel{d}{=} 
(\frac{1}{2} L_{u/2})_{u}$, 
where $L$ denotes the local time process of a normalized 
Brownian excursion. The reason is: Up to a factor 1/2 
$L$ is related to a normalized Brownian 
excursion in the same way, as $H$ is derived above from $Y$, which 
follows from a result of Jeulin \cite{je} (compare Biane 
\cite{bi}, Th\'{e}or\`{e}me 3). This version of Theorem 1 has been observed 
by Aldous in special cases (as the geometric offspring distributions) and
conjectured in the general finite variance case, compare
\cite{al2}. A first proof of the conjecture was given by 
Drmota and Gittenberger \cite{dr1}. They mastered the formidable task 
to obtain convergence of the finite-dimensional distributions as well 
as tightness, using 
generating functions and thereby generalizing work of Kennedy 
\cite{ken} on the one-dimensional distributions. Pitman 
\cite{pi} surrounded the difficulties by imbedding 
the problem into the context of convergence of strong solutions of 
stochastic differential equations. The limiting distribution of the 
maximum of the height profile, the `width' of the family tree, has 
been found in the finite variance case already by Tak\'{a}cs \cite{tk}.

The relations to stochastic analysis can be further developed in the 
case $\sigma^{2} < \infty$. Note 
that the defining equation for $C$ can be written as 
\begin{displaymath}
    u = \int_{0}^{C_{u}} \frac{dt}{Y_{t}} \mbox{ for } 
	0 < u \leq h, \mbox{ and } C_{u} = 1 \mbox{ for } u > h,
\end{displaymath}
where
\begin{displaymath}
	h = \int_{0}^{1} \frac{dt}{Y_{t}}.
\end{displaymath}
$h$ is the asymptotic height of the rescaled tree, as will become 
clear in the next section. By differentiation with respect to $u$ we get
\begin{displaymath}
    H_{u} =  \frac{dC_{u}}{du} = Y(C_{u}).
\end{displaymath}
Now a Brownian excursion $Y$ solves the stochastic equation $dY = dW + 
(1-\frac{Y^{2}}{1-t})\frac{dt}{Y}$, with a standart Brownian motion $W$ 
(compare \cite{ro}, chapter IV, section (40.4)). 
Viewing $C$ as a time-change, this leads to the stochastic equation
\begin{displaymath}
    dH = \sqrt{H} dB +\left(1-\frac{H^{2}}{1-C}\right)dt,
\end{displaymath}
with a standart Brownian motion $B$. Pitman \cite{pi} also obtains this 
equation  and discusses it in detail, 
therefore there is no need to trace this aspect further. The idea of 
using $C$ as a time-change goes back to Lamperti \cite{la}.

The mentioned proofs all focus on the relationship to Brownian local time. 
In contrast we shall not rely on local times in this paper. (A short 
description of our proof in the finite variance case already appeared in the 
technical report \cite{ker}.) The reason is that 
in the case $\alpha < 2$ the connection to local times breaks 
down. Then $Y$ and consequently $H$ exhibits a.s. jumps. Since
local time processes of L\'{e}vy-processes (if existent) have a.s. continuous 
paths (compare f.e. \cite{ber}, chapter V.1), they are no longer 
suited.

This can be explained on a heuristic level, too. Consider the 
following construction, going back to Harris \cite{ha} and used by 
different people. Traverse 
the individuals of $T$ in the following manner: From individual 
$i$  pass over to its oldest child, which has not yet been visited, 
resp. return to its predecessor, if all children of $i$ have already 
been visited. This gives a traversal through $T$ with the root as starting 
and end point. Each edge is passed twice, once forward and once backward. Next 
consider the associated random path, which increases one unit, if we change over 
to a child, and decreases one unit, if we go back to a predecessor. 
In the case of a geometric offspring distribution we get a true 
random walk excursion, conditioned to return to zero after $2n$ steps for a 
$CGW(n)$-tree. This is due to the lack of memory of the geometric 
distribution. The number of upcrossings from level $k-1$ to level 
$k$ is equal to the size $Z_{k}$ of the $k$'th generation, which 
makes the relation to local times obvious in this situation. In general 
the random path exhibits complicated dependence properties. In 
the case $\sigma^{2} < \infty$ they are of a local nature and vanish 
in the limit $n \rightarrow \infty$, as was shown by Aldous, 
such that the asymptotic height 
profile can still be described by Brownian local time. For 
$\alpha < 2$ however, the dependence structure survives
in the limit.

The combinatorics of $CGW(n)$-trees have been widely studied by means 
of generating functions (see f.e. \cite{dr1,fl,ken,me}). Probabilistic 
methods have been introduced in Kolchin \cite{ko} and in particular 
by Aldous \cite{al1,al2}. Our proof of theorem 1 is based on two 
probabilistic constructions, which are valid for the infinite variance 
case, too. The first one will be described in section 2, it
establishes a connection between Galton-Watson trees and 
suitable random walk excursions. Though the relationship has been known for 
quite a while (compare \cite{dw}), the scope of this approach 
has been enlarged considerably only recently (see \cite{ben,bo,ga}). 
In our context it allows to reduce  convergence of $H^{n}$ to convergence of 
excursions. The required continuity theorem will be 
developed in section 3. The second probabilistic construction, 
which will be presented in section 4, is size-biasing of Galton-Watson trees. 
We use it to describe the bottom of the trees, which is of some 
interest of its own and will help 
to check a main condition of the continuity theorem. This 
concept goes back to Geiger \cite{ge}, who developed a construction 
due to Lyons, Pemantle and Peres \cite{ly}. Section 5 addresses the question 
of convergence of excursions. 

Thus the height profile will be considered as a functional of a 
random walk excursion $S$. Other quantities of the tree can be viewed as 
well as functionals of $S$. In this manner we can also treat 
the height profile of random forests, as discussed by Drmota and 
Gittenberger \cite{dr2} and Pitman \cite{pi} in the finite variance 
case. We contend ourselves by stating a version of the theorem, which 
is valid for infinite variances, too. A \emph{conditioned random 
forest} consists of $l$ Galton-Watson trees, conditioned to contain 
altogether $n$ individuals. Let now $Z_{k}$ be the total number of 
all individuals in generation $k$ in one of the $l$ trees. Suppose 
that $l \sim \gamma  a_{n}$, as $n \rightarrow \infty$, with $\gamma 
\geq 0$. Then the height profile $(H_{u}^{n})_{u} = 
(a_{n}^{-1}Z_{[nu/a_{n}]})_{u}$ converges in distribution. The limiting 
process can be described as follows. Let $X_{\gamma}=(X_{\gamma,s})_{0\leq 
s \leq 1}$ be a L\'{e}vy-process as above, now conditioned to hit $-\gamma$ 
at the moment $s=1$ for the first time. It is built up from the excursions
 \begin{displaymath}
 	Y_{\eta} = (Y_{\eta,s})_{s \leq L_{\eta}} = 
(X_{s}+\eta)_{T_{\eta} \leq s \leq T_{\eta+}},\; 0 \leq \eta 
\leq \gamma,
 \end{displaymath} 
with $T_{\eta} = \inf\{ s: \: X_{s} = - \eta\}$, 
$T_{\eta+} = \inf\{s : \: X_{s} < - \eta\}$ and $L_{\eta} = T_{\eta+}- 
T_{\eta}$. Then the limiting process is given by
\begin{displaymath}
	\sum_{0 \leq \eta \leq \gamma} H_{\eta,u}, \; u\geq 0,
\end{displaymath}
where $H_{\eta,u}$ is derived similarly as above from
\begin{displaymath}
	C_{\eta,u} = \sup \{s \leq L_{\eta} \: : \: \int_{0}^{s} 
	\frac{dt}{Y_{\eta,t}}\leq u \}.
\end{displaymath}
The sum is a.s. finite for every $u>0$, which reflects the fact, that 
also in the limit only finitely many trees contribute to the height 
profile. This result can be 
proved in much the same manner (and with only little additional 
effort), as we shall obtain Theorem 1 below. For $\gamma = 0$ we are 
back in the situation of Theorem 1.

\section{Trees and Random Walk Excursions}

It has been known for some time, that Galton-Watson trees can be 
imbedded into random walks (compare \cite{bo} and the references therein). 
This is implicit in Dwass' important paper \cite{dw}.
Here  we give a combinatorical treatment, which allows generalization.
Let $T$ be a tree of size $n$. Suppose that we label 
the individuals  in $T$ with the numbers $1,2, \ldots ,n$ 
such that the root gets label 1, the individuals in the first 
generation the labels $2,\ldots,Z_{0}+Z_{1}$ (say from left to right), the 
individuals in the second generation the labels $Z_{0}+Z_{1}+1, \ldots, 
Z_{0}+Z_{1}+Z_{2}$ and so forth. Using these labels we define a random path 
$S=(S(0), S(1), \ldots, S(n))$ recursively by 
\begin{displaymath}
	S(0)=1, \; S(i) = S(i-1) + \xi_{i} - 1, \; i= 1, \ldots n,
\end{displaymath}
where $\xi_{i}$ denotes the offspring number of the individual with label $i$.
We can imagine that the path arises as follows: To its $i$'th 
increment individual $i$ contributes the downward step $-1$, whereas 
each of its children contributes one step $+1$ upwards. Then each 
individual is responsible for one upward and one downward step, except 
the root, which has no predecessor and thus contributes only a step 
downwards. Since $S(0) = 1$,
\begin{displaymath}
	S(n) =0. 
\end{displaymath}
Furthermore, individuals always have smaller labels than their 
offspring, therefore the upward step of an individual appears 
before its downward step. Clearly this implies 
\begin{displaymath}
	S(i) > 0 \mbox{ for all } i < n,
\end{displaymath}
i.e. $S$ is an \emph{excursion} of length $n$. Conversely, given such
an excursion $S$ of length n, with $S(i) \geq S(i-1) - 1$, we 
can construct a tree, fitting to the excursion $S$. Namely, from $S$ 
we read off $\xi_{i} = S(i) - S(i-1) + 1$, and the given labelling rule allows 
us to grow the tree from its root. Thus there is a one-to-one 
correspondence between trees of size $n$ and excursion of length 
$n$.

As to the probabilistic aspect of the construction note, that for a 
Galton-Watson tree the $\xi_{i}$ are independent random 
variables, such that $S$ becomes an ordinary random walk excursion. 
Likewise $S$ is a random walk excursion of duration n, if $T$ is a 
$CGW(n)$-tree. \bigskip

\textbf{Remarks} \quad 1.) These considerations remain valid for other 
ways of labelling. One possibility is to label according to 
depth-first search, which has been exploited in 
\cite{ben,ga}. In general the following properties are 
required:\\
i) $1$ is the label of the root.\\
ii) The individuals with labels $1, \ldots,i$ form a subtree for 
any $i < n$. In other words: Any individual has a smaller label than 
any of its children.\\
iii) Given the subtree with labels $1, \ldots, i$ and the numbers 
$\xi_{1}, \ldots, \xi_{i}$ there is a rule, which specifies, which child of the 
individuals $1, \ldots ,i$ gets the label $i+1$.\\
2.) A random forest of $l$ trees and $n$ individuals can be described by a 
random walk path $S$ with $S(0)= l$, $S(i) >0$ for $i<n$ 
and $S(n)=0$. If we label the trees one after the other, then the 
$r$'th tree is represented by the part of the random walk between the 
hitting times of $l+1-r$ and $l-r$. \bigskip

In our labelling $1, \ldots , Z_{0}+ Z_{1}+ \ldots + Z_{k-1}$ are 
just the members of generation 0 to $k-1$. They contribute a negative 
step to $S(Z_{0}+ Z_{1}+ \ldots + Z_{k-1})$. Their children, i.e the 
individuals in generation 1 to $k$, add a positive step. Therefore
\begin{displaymath}
	Z_{k} = S( Z_{0}+ Z_{1}\ldots + Z_{k-1}).
\end{displaymath}
This is the announced random walk representation for the height 
profile, which has been used by several authors (see \cite{bo}). 
We transform it into a differential equation. Define
\begin{displaymath}
	C_{u}^{n} = \int_{0}^{u} a_{n}^{-1} \: Z_{[nv/a_{n}]} \; dv,
\end{displaymath}
in particular
\begin{displaymath}
	C_{ka_{n}/n}^{n} = \frac{1}{n} (Z_{0}+ \ldots + Z_{k-1}).
\end{displaymath}
Further let for $0 \leq s \leq 1$
\begin{displaymath}
	S^{n}(s) = a_{n}^{-1} S([ns])
\end{displaymath}
and
\begin{displaymath}
	Y^{n}(s) = S^{n}(C_{ka_{n}/n}^{n}) \mbox{ for } s \in 
	[C_{ka_{n}/n}^{n},C_{(k+1)a_{n}/n}^{n}).
\end{displaymath}
Then the above equation translates into the ordinary differential 
equation
\begin{displaymath}
	\frac{d}{du} C_{u}^{n} = Y^{n}(C_{u}^{n}),
\end{displaymath}
which by integration leads to
\begin{displaymath}
	u = \int_{0}^{C_{u}^{n}} \frac{dt}{Y^{n}(t)} \mbox{ for } u \leq h_{n} 
	\mbox{ and } C_{u}^{n} =1 \mbox{ for } u > h_{n},
\end{displaymath}
with
\begin{displaymath}
	h_{n} = \int_{0}^{1} \frac{dt}{Y^{n}(t)}.
\end{displaymath}
$h_{n}$ obviously is the height of the rescaled tree. --- It is 
now our plan to reduce the question of convergence 
of $H^{n}$ and $C^{n}$ to that of $S^{n}$ and $Y^{n}$. The next 
section provides the required continuity statement.

\section{A Continuity Theorem}

Let $D$ and $D'$ be the spaces of all c\'{a}dl\'{a}g-functions 
$f:[0,1] \rightarrow \mathbb{R}$ resp. $g:[0,\infty) \rightarrow 
\mathbb{R}$. By endowing them with the common Skorohod-distance 
(compare \cite{et}) we make them to complete metric spaces. Thus
functions $f_{n}$ converge to $f$ in $D$, if there are 
increasing bijections $\alpha_{n}: [0,1] \rightarrow [0,1]$, such 
that, as $n \rightarrow \infty$,
\begin{displaymath}
	\sup_{t} |\alpha_{n}(t) - t | \rightarrow 0, \mbox{\qquad} \sup_{t} 
	|f_{n}(t) - f(\alpha_{n}(t))| \rightarrow 0.
\end{displaymath}
Similarly $g_{n} \rightarrow g$ in $D'$, if there are 
increasing bijections $\beta_{n}: [0,\infty) \rightarrow [0,\infty)$, such 
that for all $v < \infty$
\begin{displaymath}
	\sup_{u \leq v} |\beta_{n}(u) - u | \rightarrow 0, \mbox{\qquad} 
	\sup_{u \leq v} 
	|g_{n}(u) - g(\beta_{n}(u))| \rightarrow 0.
\end{displaymath}
If $g$ is a continuous function, we may choose $\beta_{n}(u) =u$.
 
We begin by collecting some analytical facts. Let $D_{+}$ be the space 
of non-negative $f \in D$.

\begin{lemma} Let $f,f_{n} \in D_{+}$.\\
i) For given $s$ the mapping $f \mapsto \int_{0}^{s} 
dt/f(t)$ (with the possible value $\infty$)
is lower-semicontinuous (and thus measurable).\\
ii) If $\inf_{a \leq t \leq b} f(t) > 0$ for some $0 \leq a < b \leq 
1$, and if $f_{n} \rightarrow f$, then 
\begin{displaymath}
	\sup_{a \leq s \leq b} \left| \int_{a}^{s} \frac{dt}{f(t)} - 
	\int_{a}^{s} 
	\frac{dt}{f_{n}(t)} \right| \rightarrow 0.
\end{displaymath}
\end{lemma}

\emph{Proof} \quad Let $\alpha_{n}$ as above. As is well-known,
it may be assumed without loss of generality, that they 
are differentiable functions such that
\begin{displaymath}
	\sup_{t} |\alpha_{n}^{'}(t) -1 | \rightarrow 0.
\end{displaymath}
 
i) \quad If $f_{n} \rightarrow f$, then for $\delta > 0$ by Fatou's Lemma
\begin{displaymath}
	\int_{0}^{s-\delta} \frac{dt}{f(t)} \leq \liminf_{n} 
	\int_{0}^{s-\delta} \frac{dt}{f_{n}(\alpha_{n}^{-1}(t))} \leq
	\liminf_{n} \int_{0}^{s} \frac{\alpha_{n}^{'}(t) \: dt}{f_{n}(t)}
\end{displaymath}
and consequently
\begin{displaymath}
	\int_{0}^{s} \frac{dt}{f(t)} \leq \liminf_{n} \int_{0}^{s} 
	\frac{dt}{f_{n}(t)}.
\end{displaymath}
 
ii) Due to uniform convergence of the integrands
\begin{displaymath}
	\sup_{a \leq s \leq b} \left| \int_{a}^{s} 
	\frac{dt}{f_{n}(t)} - \int_{a}^{s} 
	\frac{dt}{f(\alpha_{n}(t))} \right| \rightarrow 0.
\end{displaymath}
Further, substituting $\alpha_{n}(t) = w$
\begin{displaymath}
	\int_{a}^{s} \frac{dt}{f(\alpha_{n}(t))} = 
	\int_{a}^{\alpha_{n}(s)} 
	\frac{dw}{\alpha_{n}^{'}(t)f(w)} \rightarrow \int_{a}^{s} 
	\frac{dw}{f(w)}
\end{displaymath}
uniformly for all $s \in [a,b]$. \hfill $\Box$ \bigskip

Define now for $f \in D_{+}$ a function $g = \phi(f) \in D'$ by
\begin{displaymath}
	g(u) = \sup \{s \leq 1 \: | \: \int_{0}^{s} \frac{dt}{f(t)} \leq u\}.
\end{displaymath}
Thus as above 
\begin{displaymath}
	u = \int_{0}^{g(u)} \frac{dt}{f(t)} \mbox{ for } u < h, \mbox{ and } 
	g(u) = 1 \mbox{ for } u \geq h,
\end{displaymath}
with
\begin{displaymath}
	h = \int_{0}^{1} \frac{dt}{f(t)}.
\end{displaymath}
$g$ is continuous and increasing. It is everywhere differentiable from the right 
(since $f$ is continuous from the right), and the derivative is given by 
\begin{displaymath}
     \frac{d^{+}}{du} g(u) = f(g(u)) \mbox{ for } u < h \mbox{ 
     and } \frac{d^{+}}{du} g(u) = 0 \mbox{ for } u \geq h.
\end{displaymath}
We denote
\begin{displaymath}
	\psi(f) = \frac{d^{+}}{du} \phi(f).
\end{displaymath}

\begin{lemma}
Suppose $f_{n} \rightarrow f$ in $D_{+} $, $ \sup_{u} |\phi(f_{n})(u) - 
\phi(f)(u)| \rightarrow 0$ and $\int_{0}^{s} dt/f(t) < \infty$ for 
all $s<1$. Then $\psi(f_{n}) \rightarrow \psi(f)$ in $D'$.
\end{lemma}

\emph{Proof} \quad  Denote $g = \phi(f)$, $g_{n} = \phi(f_{n})$ and 
$h_{n} = \int_{0}^{1} dt/f_{n}(t)$.
Let $\alpha_{n}(t)$ be as above. The required bijections $\beta_{n}: 
\mathbb{R}_{+} \rightarrow \mathbb{R}_{+}$ are defined as
\begin{displaymath}
	\beta_{n}(u) = \left\{
	\begin{array}{ll}
		g^{-1}(\alpha_{n}(g_{n}(u))) & \mbox{ for } u \leq h_{n},  \\
		\beta_{n}(h_{n}) + (u-h_{n}) & \mbox{ for } u > h_{n}.
	\end{array}
	\right.
\end{displaymath}
Then
\begin{eqnarray*}
	\sup_{u \leq h_{n}} |g(\beta_{n}(u)) - g_{n}(u)| & = & \sup_{u \leq 
	     h_{n}} | \alpha_{n}(g_{n}(u)) - g_{n}(u)|  \\
	 & \leq & \sup_{t} |\alpha_{n}(t) - t| \rightarrow 0.
\end{eqnarray*}
By assumption it follows
\begin{displaymath}
	\sup_{u \leq h_{n}} |g(\beta_{n}(u)) - g(u)| \rightarrow 0.
\end{displaymath}
If $h = \int_{0}^{1} dt/f(t) < \infty$, then $g^{-1}(s) = \int_{0}^{s} 
dt/f(t)$ is uniformly continuous on $[0,1]$, and it follows
\begin{displaymath}
	\sup_{u} |\beta_{n}(u) - u| = \sup_{u \leq h_{n}} |\beta_{n}(u) - u| 
	\rightarrow 0.
\end{displaymath}
If on the other hand $h = \infty$, then $h_{n} 
\rightarrow \infty$ and $g^{-1}$ is uniformly continuous on every 
intervall $[0,s]$ with $s<1$. In this case we may conclude
\begin{displaymath}
	\sup_{u \leq v} |\beta_{n}(u) - u| \rightarrow 0
\end{displaymath}
for every $v>0$. This is one of the desired properties.

Next for $u>h_{n}$ we have $1=g_{n}(h_{n})=g_{n}(u)$ and $1= 
\alpha_{n}(1) = \alpha_{n}(g_{n}(h_{n})) = g(\beta_{n}(h_{n})) 
=g(\beta_{n}(u))$, therefore
\begin{eqnarray*}
	\lefteqn{\sup_{u} |f_{n}(g_{n}(u)) - f(g(\beta_{n}(u)))|} \\
	 & = & \sup_{u \leq h_{n} } |f_{n}(g_{n}(u)) - f(g(\beta_{n}(u)))| \\
	 & = & \sup_{u \leq h_{n}} |f_{n}(g_{n}(u)) - f(\alpha_{n}(g_{n}(u)))| \\
	 & = & \sup_{t} |f_{n}(t) - f(\alpha_{n}(t))| \rightarrow 0.
\end{eqnarray*}
Since $\psi(f) = f \circ g$ and $\psi(f_{n}) = f_{n} \circ g_{n}$, 
also
\begin{displaymath}
	\sup_{u} |\psi(f_{n})(u) - \psi(f)(\beta_{n}(u))| \rightarrow 0,
\end{displaymath}
which proves the claim.  \hfill $\Box$ \bigskip

We are now ready to prove the main result of the section. Let 
$Y,Y^{n}$  be stochastic processes with paths in $D_{+}$, and define 
processes $C, C^{n}, H$ and $H^{n}$ by
\begin{displaymath}
\begin{array}{ll}
	C_{u} = \phi(Y)(u), &   C_{u}^{n} = \phi(Y^{n})(u), \\
	H_{u} = \psi(Y)(u), &  H_{u}^{n} = \psi(Y^{n})(u).
\end{array}
\end{displaymath}

\begin{theorem} 
Assume $\inf_{\delta \leq s \leq 1- \delta} Y_{s} > 0$
a.s. for all $\delta > 0$ and, as $\epsilon  \rightarrow 0$,
\begin{displaymath}
	\limsup_{n} \mathrm{P}^{n}(C_{u}^{n} \leq \epsilon) \rightarrow 0
\end{displaymath}
for all $u>0$. Then $C_{u} > 0$ a.s. for all $u > 0$, i.e. $H$ is 
non-degenerate, and $C^{n} \rightarrow 
C$, $H^{n} \rightarrow H$ in distribution, as $n \rightarrow 
\infty$.
\end{theorem}

\emph{Proof} \quad Due to a wellknown theorem of Skorokhod (see 
\cite{po}, chapter IV.3, Theorem 13) we may 
assume that the processes $Y_{n}$ and $Y$ are defined on a single 
probability space $(\Omega, \cal{A}, \mathrm{P})$, and that $Y_{n} 
\rightarrow Y$ a.s. in $D_{+}$. Because of semicontinuity (Lemma 2 i)) 
and the definition of $C$ it follows $C_{u} \geq \limsup_{n} C_{u}^{n}$ a.s.. 
From Fatou's Lemma  
\begin{displaymath}
	\mathrm{P}(C_{u} < \epsilon ) \leq \mathrm{P}(\limsup_{n} C_{u}^{n} 
	< \epsilon) \leq \limsup_{n} \mathrm{P}( C_{u}^{n} < \epsilon). 
\end{displaymath}
Thus our assumptions imply $C_{u} > 0$ a.s. for all $u > 0$.

Next let $C_{u} < 1$. Then $C_{u}^{n} < 1$ for large $n$, therefore $u= 
\int_{0}^{C_{u}} dt/Y_{t} = \int_{0}^{C_{u}^{n}} dt/Y_{t}^{n}$. It 
follows
\begin{eqnarray*}
    |C_{u}^{n} - C_{u}| & \leq  &  \max Y^{n} \; \left| 
	\int_{0}^{C_{u}^{n}} \frac{dt}{Y_{t}^{n}} - \int_{0}^{C_{u}} 
	\frac{dt}{Y_{t}^{n}} \right| \\
	 & = & \max Y^{n} \; \left| \int_{0}^{C_{u}} \frac{dt}{Y_{t}} 
	 - \int_{0}^{C_{u}} \frac{dt}{Y_{t}^{n}} \right|.
\end{eqnarray*}
Therefore for any $a<1$ and $b$ such that $C_{b} = a$
\begin{eqnarray*}
	\sup_{u}|C_{u}^{n}-C_{u}| & \leq &  \sup_{u \leq b} 
	|C_{u}^{n}-C_{u}| + (1-C_{b}^{n}) + (1-C_{b})  \\
	 & \leq & 2 \sup_{u \leq b} |C_{u}^{n}-C_{u}| + 2(1-a)  \\
	 & \leq & 2 \max Y^{n} \; \sup_{s \leq a } \left| \int_{0}^{s} 
	 \frac{dt}{Y_{t}} - \int_{0}^{s} \frac{dt}{Y_{t}^{n}} \right|+ 2(1-a).
\end{eqnarray*}
Now $\max Y^{n} \rightarrow \max Y$ a.s.. In view of Lemma 2 ii) 
it follows for every $ \eta > 0$ there is a $v >0$ 
(chose $a$ sufficiently close to 1), such that for all $\epsilon > 0$
\begin{eqnarray*}
	\lefteqn{\limsup_{n} \; \mathrm{P}(\sup_{u} |C_{u}^{n} - C_{u}| > \eta)} \\
	 & \leq & \eta + \limsup_{n} \; \mathrm{P}(\left|\int_{0}^{\epsilon}  
	 \frac{dt}{Y_{t}} - \int_{0}^{\epsilon} \frac{dt}{Y_{t}^{n}} 
	 \right| > v)  \\
	   & \leq & \eta + \mathrm{P}(C_{v} < \epsilon) + \limsup_{n} \; 
	   \mathrm{P}(C_{v}^{n} < \epsilon)   \\
	 & \leq & \eta + 2\; \limsup_{n} \;  P(C_{v}^{n} < \epsilon).
\end{eqnarray*}
By assumption the righthand side can be made arbitrarily small, such 
that $\sup_{u} |C_{u} - C_{u}^{n}| \rightarrow 0$ in probability. In 
view of Lemma 3 $H^{n} \rightarrow H$ in probability, 
and the proof is finished.
\hfill $\Box$ \bigskip

\section{Size-biased Galton--Watson Trees}

In this section we verify one of the conditions of Theorem 4 for 
height profiles of a $CGW(n)$-tree. We shall make use of 
size-biased trees, as they have been constructed by Geiger for binary 
branching trees resp. splitting trees (compare \cite{ge}). Let
\begin{displaymath}
	\varphi(\lambda) = \sum_{x=0}^{\infty} \lambda^{x} p_{x}
\end{displaymath}
be the generating function of the offspring distribution $(p_{x})$. Fix 
$\lambda \in (0,1)$ and consider the probability weights
\begin{displaymath}
	q_{x} = \lambda^{x}p_{x}/\varphi(\lambda).
\end{displaymath}
Then we have for a Galton Watson tree $T$ the formula
\begin{displaymath}
	\mathrm{Q}(T=t) = \prod_{i} q_{d(i)} = \lambda^{s(t)-1} \varphi(\lambda)^{-s(t)} 
	\prod_{i}p_{d(i)} =  \lambda^{s(t)-1} \varphi(\lambda)^{-s(t)}\: 
	\mathrm{P}(T=t),
\end{displaymath}
where $d(i)$ denotes the number of children of individual $i$ in the 
(nonrandom) tree $t$, and $\mathrm{P}$ and $\mathrm{Q}$ denote the probability 
measures, corresponding to the offspring distributions $(p_{x})$ 
resp. $(q_{x})$. It follows the known fact (compare \cite{ken}) that two 
$CGW(n)$-tree with offspring distributions $(p_{x})$ and $(q_{x})$ 
are equal in distribution. Passing over to $(q_{x})$ has the effect 
that the total size of $T$ becomes finite in mean 
(\cite{fe1}, chapter XII.5):
\begin{displaymath}
	\mathrm{E}_{\mathrm{Q}}s(T) = (1-\mu)^{-1} < \infty,
\end{displaymath}
where
\begin{displaymath}
	\mu = \sum_{x} xq_{x} < 1.
\end{displaymath}
Thus (differently from $\mathrm{P}$) we may bias $\mathrm{Q}$ by 
introducing the probability measure
\begin{displaymath}
	\widehat{\mathrm{Q}}(T=t) = (1-\mu) s(t) \: \mathrm{Q}(T=t).
\end{displaymath}
Note: Conditioning on size $s(T)=n$, it makes no difference, whether we 
consider $\mathrm{Q}$ or $\widehat{\mathrm{Q}}$. If moreover $M$ is chosen 
purely at random from the vertices of $T$, i.e. from $\{1,2,\ldots, 
s(T)\}$, then
\begin{displaymath}
	\widehat{\mathrm{Q}}(T=t,M=m) = (1-\mu) \mathrm{Q}(T=t).
\end{displaymath}
Let $(\widehat{q}_{x})$ be the distribution obtained from $(q_{x})$ by 
size-biasing, i.e. 
\begin{displaymath}
	\widehat{q}_{x} = x q_{x}/\mu,
\end{displaymath}
then we may rewrite the above formula as 
\begin{displaymath}
	\widehat{\mathrm{Q}}(T=t,M=m) = (1-\mu)\mu^{g} 
	\prod_{i<m} \widehat{q}_{d(i)} 
	\frac{1}{d(i)} \; \prod_{i \not< m} q_{d(i)},
\end{displaymath}
where $g$ denotes the generation of individual $m$, and $i<m$ means 
that individual $i$ is a predecessor of $m$. Following Geiger 
\cite{ge}, this formula gives raise to a probabilistic construction 
for the size-biased Galton Watson tree: \bigskip  

\textbf {Construction of the tree $\widehat{T}$}
\begin{itemize}
\item Let $G$ be a random variable with geometric distribution and 
parameter $\mu$. Choose independent random variables  
$\widehat{\xi}_{1},\ldots, \widehat{\xi}_{G}$  with 
distribution $(\widehat{q}_{x})$. Let $\zeta_{j}$ be random numbers, 
taken independently and uniformly from $\{1,2,\ldots,
\widehat{\xi}_{j}\}$, $j=1,\ldots,G$.\\
\item  The trunk of $\widehat{T}$: Build up a line of $G$ consecutive 
individuals, where the $j$'th individual has 
$\widehat{\xi}_{j}$ children. Let the first one be the 
tree's root and the $(j+1)$'th one be the $\zeta_{j}$'th child 
(in the order of birth)  of the j'th individual. Give the label $M$ to the 
$\zeta_{G}$'th child of the $G$'th individual.\\
\item  The tree top of $\widehat{T}$: Besides these $G$ individuals 
the trunk contains $y=\widehat{\xi}_{1}+\ldots+\widehat{\xi}_{G} +1-G$ 
additional individuals. They 
propagate in the usual Galton-Watson manner. This means: In order to 
complete the tree $\widehat{T}$ we attach independent Galton-Watson trees 
$T_{1},\ldots,T_{y}$ with offspring 
distribution $(q_{x})$ to the trunk.
\end{itemize}

In this manner we obviously obtain a size-biased tree together with 
an individual $M$, taken at random from the tree, and belonging to 
generation $G$:
\begin{displaymath}
	\mathrm{Q}(\widehat{T}=t, M=m) = \widehat{\mathrm{Q}}(T=t,M=m).
\end{displaymath}
Note that for this construction it is only required 
that $(p_{x})$ has mean 1. In the limiting case $\lambda =1$  (or 
equivalently $\mu = 1$) $G$ takes a.s. the value $\infty$, and 
$\widehat{q}_{x}$ is equal to
\begin{displaymath}
	\widehat{p}_{x} = xp_{x}.
\end{displaymath}
Then we get an \emph{infinite size-biased tree} $\widetilde{T}$, which 
already appeared in the work of Grimmett \cite{gr}, Kesten 
\cite{kes}, Aldous \cite{al1} 
and others. We use this tree to describe the asymptotic shape of the 
lower part of a $CGW(n)$-tree. For any tree $t$ let $t(k)$ be the tree, 
which results by cutting off all individuals in $t$ belonging to a 
generation greater than $k$. Thus, $t(k)$ has in generation $k'$ the 
generation sizes $z_{k'}$ for $k' \leq k$, and $0$ for $k' > k$. 
Further let $\mathcal{T} (k)$ be the set of trees with height at most $k$.
\begin{theorem}
Let $T$ be a Galton-Watson tree with an offspring distribution $(p_{x})$
fulfilling assumption A, and let $(k_{n})$ be a sequence of 
natural numbers such that $k_{n} =o(n/a_{n})$. Then for $n \rightarrow \infty$
\begin{displaymath}
	\sup_{B \subset  \mathcal{T} (k_{n})} \: | \mathrm{P}(T(k_{n}) \in B \:
	| \: s(T) =n) - \mathrm{P}(\widetilde{T}(k_{n}) \in B) | \rightarrow 0.
\end{displaymath}
\end{theorem}
In the special case of a Poisson offspring distribution this was proved 
by Aldous in \cite{al1}. We prepare the proof by gathering some facts 
of an analytical character.  Let
\begin{displaymath}
	v(x) = \sum_{y \leq x} y(y-1)p_{y}.
\end{displaymath}
In order that $(p_{x})$ belongs to the domain of attraction of a 
stable law with index $\alpha$, it is necessary and 
sufficient, that $v(x)$ varies regularly at infinity with exponent 
$2-\alpha$. Then $a_{n}$ can be chosen as any sequence with the 
property
\begin{displaymath}
    a_{n}^{2} / v(a_{n}) = cn (1+o(1))
\end{displaymath}
with some $c>0$ (compare \cite{fe2}, chapter XVII.5). 
It follows that $a_{n}$ varies 
regularly with exponent $1/ \alpha$. We have a local limit law at 
our disposal, which in our case reads as follows.
\begin{proposition}
Let $\xi_{1}, \xi_{2},\ldots$ be i.i.d. random variables with common 
distribution $(p_{x})$, satisfying assumption A. Then uniformly in $x \in 
\mathbb{Z}$
\begin{displaymath}
	\mathrm{P}(\xi_{1}+\ldots+\xi_{n}-n= x) = 
	a_{n}^{-1}g(x/a_{n})(1+o(1)),
\end{displaymath}
where $g$ denotes the (continuous and strictly positive) density of 
the limit law $\nu$.
\end{proposition}

The proof can be found in \cite{gn}, chapter 9, up to an exceptional 
case, and in \cite{mu} in full generality. 

\begin{lemma}
Let $\widehat{\xi}_{1}, \widehat{\xi}_{2},\ldots$ be independent random 
variables with distribution $(\widehat{p}_{x})$. Then assumption A implies
convergence of $a_{n}^{-1}(\widehat{\xi}_{1}+\ldots+ \widehat{\xi}_{[n/a_{n}]})$  
in distribution.
\end{lemma}

\emph{Proof} \quad   If $v(\infty) = \sigma^{2} < \infty$, then 
$a_{n}$ is asymptotically proportional to $n^{1/2}$ and
$(\widehat{p}_{x})$ has finite mean, and the claim follows from the 
ordinary law of large numbers. 
The case $\alpha = 2, v(\infty) = \infty$ 
is similar. Then $v(x) = \sum_{y \leq x} (y-1) \widehat{p}_{y}$ is a slowly 
varying function.  It follows $x \mathrm{P}(\widehat{\xi} \geq x) = o(v(x))$, 
as $x \rightarrow \infty$ (compare \cite{fe2}, chapter VIII.9, Theorem 2). 
This allows to apply a generalized law of large numbers, as given f.e. in 
\cite{du}, chapter 1.5, exercise 5.11.

If $\alpha < 2$, then $x \sum_{y \geq x} \widehat{p}_{y} \sim 
(2-\alpha)(\alpha-1)^{-1}v(x)$, as follows from the cited theorem in 
\cite{fe2}. Therefore $\mathrm{P}(\widehat{\xi} \geq x)$ is regular 
varying with exponent $1-\alpha$, such that $\widehat{\xi}$ belongs in 
distribution to the 
domain of attraction of a positive stable law. Moreover $\mathrm{P}
(\widehat{\xi} \geq a_{n}) \sim c a_{n}/n$ for some $c>0$, from which 
the claim follows by standard results on convergence in distribution 
to stable laws. \hfill $\Box$ \bigskip

\begin{lemma}
Let $T, T_{1},T_{2},\ldots$ be independent Galton-Watson trees with 
offspring distribution $(p_{k})$, satisfying assumption A. Then it follows:\\
i) $\mathrm{P}(s(T) = n) = n^{-1} a_{n}^{-1}g(0)(1+o(1))$, as $n \rightarrow 
\infty$. \\
ii) $n^{-1}(s(T_{1})+\ldots+s(T_{[a_{n}]}))$ converges in distribution.
\end{lemma}

\emph{Proof} i) is an immediate consequence of Proposition 6 and the 
following classical formula (see Dwass \cite{dw} and 
Kolchin \cite{ko}, Lemma 2.1.3):
\begin{displaymath}
	\mathrm{P}(s(T)= n) = \frac{1}{n} \mathrm{P} (\xi_{1}+\ldots+\xi_{n} = n-1).
\end{displaymath}
Since $a_{n}$ is regularly varying with exponent $1/\alpha$, it 
follows $\mathrm{P}(s(T) \geq n) \sim (1-\alpha^{-1})g(0) a_{n}^{-1}$. 
Consequently $s(T)$ belongs in distribution to the domain of 
attraction of a positive stable law, and ii) follows from standard 
results. 
\hfill $\Box$ \bigskip

We come to the \emph{Proof of Theorem 5}.  Together with $T$ we 
consider an individual $M$, 
chosen at random from $T$. We already pointed out that
\begin{displaymath}
	\mathrm{P} (T=t,M=m \; | \;  s(T) = n) = \widehat{\mathrm{Q}} 
	(T=t,M=m  \; | \;  s(T) = n).
\end{displaymath}
Therefore instead of $T$ we may consider the 
size-biased tree $\widehat{T}$, as obtained in the above construction. 
Let $x_{1},\ldots, x_{k}$ be offspring numbers of the lower $k$ 
individuals in the trunk (including the root). Let the $(j+1)$'th 
individual in the trunk be the $d_{j}$'th child of the $j$'th 
individual, with $d_{j} \leq x_{j}$. Let 
$t_{1},\ldots,t_{z}$, $z=x_{1}+\ldots+x_{k}-k$ be the trees, growing out 
of these offspring to the right and left of the trunk. Because of the lack 
of memory of a geometric distribution
\begin{eqnarray*}
	\lefteqn{\mathrm{Q}( \widehat{\xi}_{1}=x_{1}, \ldots, 
	\widehat{\xi}_{k}= x_{k}, } \\
	 &   & \zeta_{1}= d_{1},\ldots, 
	       \zeta_{k}=d_{k},T_{1}=t_{1}, \ldots, T_{z}=t_{z}, 
	       s(\widehat{T}) = n)  \\
	 & = & \mathrm{Q}(G \geq k) \: \widehat{q}_{x_{1}}\ldots 
	        \widehat{q}_{x_{k}} x_{1}^{-1}\ldots x_{k}^{-1} \\
	 &   & \times \: \mathrm{Q}(T_{1}=t_{1})\ldots 
	   \mathrm{Q}(T_{z}=t_{z}) \: \mathrm{Q}( s(\widehat{T}) = n-k-s) \\
	 & = & \mu^{k} \widehat{q}_{x_{1}}\ldots \widehat{q}_{x_{k}}  
	     x_{1}^{-1}\ldots x_{k}^{-1}\: 
	        \mathrm{Q}(T_{1}=t_{1})\ldots \mathrm{Q}(T_{z}=t_{z})   \\
	 &   & \times \: (1-\mu)(n-k-s) \mathrm{Q}( s(T) = n-k-s),
\end{eqnarray*}
with  $s = s(t_{1}) +\ldots+s(t_{z})$. As to the dependence on $\lambda$ 
it is easy to check, that the 
probabilities on the righthand side contain altogether the factor  
$\mu^{-k} \lambda^{n-1} \varphi(\lambda)^{n}$. Therefore there are numbers 
$c(n,\lambda)$ such that
\begin{eqnarray*}
	\lefteqn{\mathrm{Q}(\widehat{\xi}_{1}=x_{1}, \ldots, \widehat{\xi}_{k}=
	     x_{k},}\\
	 &   &\zeta_{1}=d_{1}, \ldots, \zeta_{k}=d_{k}, 
	       T_{1}=t_{1}, \ldots, T_{z}=t_{z}, s(\widehat{T}) = n)    \\
	 & = & c(n,\lambda) \: \widehat{p}_{x_{1}}\ldots 
	       \widehat{p}_{x_{k}}x_{1}^{-1} \ldots x_{k}^{-1} \: 
	      \mathrm{P} (T_{1}=t_{1})\ldots \mathrm{P}(T_{z}=t_{z})  \\
	 &   & \times (n-k-s) \mathrm{P}(s(T) = n-k-s),
\end{eqnarray*}
in particular
\begin{displaymath}
	\mathrm{Q}(s(\widehat{T}) =n) = c(n,\lambda) \: n \: \mathrm{P}(s(T) = n).
\end{displaymath}
From Lemma 8 i) we get
\begin{eqnarray*}
	\lefteqn{\mathrm{P}(\widehat{\xi}_{1}=x_{1}, \ldots, \widehat{\xi}_{k}=
	x_{k},}\\
	 &     & \zeta_{1}=d_{1},\ldots, \zeta_{k}=d_{k},
	         T_{1}=t_{1}, \ldots, T_{z}=t_{z} \; | \; s(T) =n) \\
	 & \sim & \widehat{p}_{x_{1}}\ldots \widehat{p}_{x_{k}} 
	         x_{1}^{-1}\ldots x_{k}^{-1} \: \mathrm{P}
	         (T_{1}=t_{1})\ldots \mathrm{P}(T_{z}=t_{z}),
\end{eqnarray*}
as long as $k+s(t_{1})+\ldots+s(t_{z}) = o(n)$. Now it is assumed 
that $k=k_{n} =o(n/a_{n})$. From Lemma 7 it follows, that in 
probability $z =\widehat{\xi}_{1} +\ldots + \widehat{\xi}_{k} -k_{n} = o(a_{n})$. 
Therefore in view of Lemma 8 $s(T_{1})+\ldots+s(T_{z}) = o(n)$. 
Thus our claim follows. \hfill $\Box$ \bigskip

We use now Theorem 5 to verify the condition of Theorem 4 for a height 
process.

\begin{lemma} Let $C^{n}$ be the cumulative height process of a 
$CGW(n)$-tree, as defined in section 2. If assumption A is 
satisfied, then as $\epsilon  \rightarrow 0$,
\begin{displaymath}
    \limsup_{n} \mathrm{P}(C_{u}^{n} \leq \epsilon) \rightarrow 0
\end{displaymath}
for all $u > 0$.
\end{lemma}

For the proof we need an estimate on the height of a Galton-Watson 
tree.

\begin{lemma}
Let $h(T) = \max \{k \: | \: Z_{k} > 0\}$ be the height of a 
Galton-Watson tree. Then assumption A implies $\mathrm{P} (h(T) > 
n/a_{n}) = c'\, a_{n}^{-1} (1+o(1))$ for some $c'>0$.
\end{lemma}

\emph{Proof} \quad Since $v(x)$ is regularly varying, by Karamata's Tauberian 
Theorem (compare \cite{fe2}, chapter XIII.5)  $\varphi''(\lambda)$ varies 
regularly at 1- with exponent $\alpha -2$, more precisely
\begin{displaymath}
	\varphi''(\lambda) = \sum_{x} x(x-1) 
	\lambda^{x-2}p_{x} \sim \Gamma(3-\alpha) \: v((1-\lambda)^{-1}).
\end{displaymath}
It follows that $\varphi'(\lambda)-1$ and $\varphi(\lambda)-\lambda$ 
vary regularly at 1-, too (compare \cite{fe2}, chapter VIII.9), and
\begin{displaymath}
	\varphi(\lambda) -\lambda \sim \frac{1}{\alpha(\alpha-1)} 
	(1-\lambda)^{2}  \varphi''(\lambda).
\end{displaymath}
This allows us to apply Lemma 2 from Slack \cite{sl}, which in our 
notation says
\begin{displaymath}
	\frac{\varphi(\mathrm{P}(h(T) \leq n)) - \mathrm{P}(h(T) \leq n)}
	{\mathrm{P}(h(T) > n)} \sim \frac{1}{(\alpha-1) n}.
\end{displaymath}
Combining these estimates and replacing $n$ by $n/a_{n}$ we get 
 \begin{displaymath}
 	\mathrm{P}(h(T) > n/a_{n}) \: v(\mathrm{P}(h(T) > n/a_{n})^{-1}) \sim 
 	\frac{\alpha \: a_{n}}{\Gamma(3-\alpha) \: n}.
 \end{displaymath}
Comparing this with $v(a_{n})/a_{n} \sim ca_{n}/n$, our claim follows. 
\hfill $\Box$ \bigskip

\emph{Proof of Lemma 9} \quad  We have to show that there are 
sufficiently many individuals in the bottom of the trees under 
consideration. This will be done first for the infinite size-biased 
tree $\widetilde{T}$. Let $T_{1},\ldots,T_{z}$ be those Galton-Watson trees, 
which grow out of the trunk of $\widetilde{T}$ at an individual, 
belonging to a generation less than $n/a_{n}$. Then, given $\eta > 
0$, in view of Lemma 7 there are numbers $c_{1},c_{2}$, such that 
$c_{1}a_{n} \leq z \leq c_{2}a_{n}$ with probability at least $1-\eta/3$. 
In view of Lemma 10 the number of trees $T_{i},i \leq 
c_{2}a_{n}$ with $h(T_{i}) > n/a_{n}$ has asymptotically a 
Poisson distribution. Therefore there is a number $l$, such that
\begin{displaymath}
	\limsup_{n} \mathrm{P}(h(T_{i}) > n/a_{n} \mbox{ for at least } l \mbox{ of 
	the } i \leq z) \leq  \eta /2.
\end{displaymath}
On the other hand because of Lemma 8 the number of trees 
$T_{i}, i \leq c_{1}a_{n}$ such that $s(T_{i}) > \delta n$ is 
asymptotically Poisson distributed, with a parameter going to 
$\infty$, as $\delta$ goes to 0. Therefore there is a $\delta > 0$ 
such that
 \begin{displaymath}
 	\limsup_{n} \mathrm{P}(s(T_{i}) > \delta n \mbox{ for less than } l \mbox{ of 
	the } i \leq z) \leq  \eta /2.
 \end{displaymath} 
Altogether we conclude: For any $\eta > 0$ there is a $\delta >0$, 
such that with probability at least $1-\eta$ there is a tree $T_{i}, 
i \leq z$, such that $h(T_{i}) \leq n/a_{n}$ and $s(T_{i}) > \delta n$. Since 
$T_{i}$ stems from an individual in a generation less than 
$n/a_{n}$, we see, that for large $n$ the number of individuals in the first 
$2n/a_{n}$ generations of $\widetilde{T}$ is bigger than $\delta n$ with
probability at least $1-\eta$.

Now let $\epsilon_{n} > 0$ be a sequence, converging to 0. Choose 
$\delta>0$ and define $m=m(n)$ by 
$\epsilon_{n}n=\delta m$. Then from the definition of $C^{n}$ in 
section 2 for $u>0$, if $n$ is large enough,
\begin{eqnarray*}
	\mathrm{P}(C_{u}^{n} \leq \epsilon_{n}) 
	& = & 
	\mathrm{P}(Z_{0}+Z_{1}+ \ldots +Z_{[nu/a_{n}]} \leq \delta m \: | \: 
	s(T)=n) \\
	&\leq &
	\mathrm{P}(Z_{0}+Z_{1}+\ldots+Z_{[2m/a_{m}]} \leq \delta m \: | \: 
	s(T)=n).
\end{eqnarray*}
Since $m=o(n)$, we may in view of Theorem 
5 switch over from $T$ to $\widetilde{T}$, therefore
\begin{displaymath}
	\limsup_{n} \mathrm{P}(C_{u}^{n} \leq \epsilon_{n}) \leq \eta 
\end{displaymath}
for all $\eta >0$. This holds for any sequence $\epsilon_{n}$, thus, 
as is not difficult 
to see, the claim of the lemma  follows. \hfill $\Box$ \bigskip

\section{Convergence of Excursions}

In this section we prove convergence in distribution of the rescaled random 
walk $S$, conditioned on the event $S(n)=0, S(i) > 0$ for $i<n$. 
The usual way to define a random walk excursion is, to condition the 
random walk on the event $S(n) \leq 0, S(i) >0$ for $i<n$. Since we 
deal with a random walk, skipfree to the left (i.e. steps to the left 
cannot be bigger than 1), this makes no difference. 
Convergence of normalized random walk excursions 
seems to be studied only in the finite variance case, see Kaigh \cite{ka}. 
In this section we derive the required generalization. Let as in section 2
\begin{displaymath}
	S^{n}(s) = a_{n}^{-1}S([ns]), \; 0 \leq s \leq 1.
\end{displaymath}

\begin{theorem}
Under assumption A the processes $S^{n}$, conditioned on the event 
$S(n)=0, S(i) > 0$ for $i<n$, converge in distribution to the normalized 
excursion $Y$ 
of the corresponding L\'{e}vy-process. It holds $Y_{0}=Y_{1}=0$ and 
$\inf_{\delta \leq s \leq 1-\delta} Y_{s} >0$ a.s. for all $\delta >0$. 
\end{theorem}

Before proving this result let us first complete the proof of our main 
theorem.\bigskip 

\emph{Proof of Theorem 1} \quad Recall the notation introduced in 
section 2. In view of Theorem 4 and Lemma 9 it remains to show 
that also the processes $Y^{n}$ converge 
in distribution. From the random walk representation of section 2 it 
follows
\begin{displaymath}
	\sup_{k} |C_{(k+1)a_{n}/n}^{n} - C_{ka_{n}/n}^{n}| = \frac{1}{n} 
	\sup_{k} Z_{k} \leq \frac{a_{n}}{n} \sup_{t} S^{n}(t).
\end{displaymath}
Because of Theorem 11 the righthand supremum converges in distribution. 
Since $a_{n} = o(n)$, we obtain
\begin{displaymath}
	\sup_{k} |t_{k+1}^{n}-t_{k}^{n}| = o(1)
\end{displaymath}
in probability, with $t_{k}^{n} = C_{ka_{n}/n}^{n}$. Without loss we may 
again assume that this convergence as well as the convergence of 
$S^{n}$ to $Y$ takes place in the a.s. sense. This means, that there 
are a.s. functions $\alpha_{n}:[0,1] \rightarrow [0,1]$ such that
\begin{displaymath}
	\sup_{t} |S^{n}(t) - Y(\alpha_{n}(t))| \rightarrow 0, \; \sup_{t} 
	|\alpha_{n}(t)-t| \rightarrow 0.
\end{displaymath}
The role of the $\alpha_{n}$ is, as is well-known, to match the jump points 
$s_{1}, s_{2}, \ldots$ of $Y$ to those jumps of $S^{n}$, which 
asymptotically are not negligible. These are given by 
$\alpha_{n}(s_{j})$. Define 
$\beta_{n}(s_{j}) = t_{k+1}^{n}$, if $\alpha_{n}(s_{j}) \in 
[t_{k}^{n},t_{k+1}^{n})$ and $j \leq d_{n}$. If $d_{n}$ is 
going sufficiently slowly to $\infty$, we obtain by linear interpolation
bijections $\beta_{n}: [0,1] \rightarrow [0,1]$, which match the jump 
points of $Y$ to those of $Y^{n}$. This implies, as is not difficult 
to see, that a.s.
\begin{displaymath}
	\sup_{t} |Y^{n}(t) - Y(\beta_{n}(t))| \rightarrow 0, \;  \sup_{t} 
	|\beta_{n}(t)-t| \rightarrow 0.
\end{displaymath}
Thus $Y^{n}\rightarrow Y$ a.s. in the Skorohod sense, and our 
claim follows. \hfill $\Box$ \bigskip

The plan of the proof of Theorem 11 is, to reduce the theorem to 
convergence of random walks and random walk bridges. The first 
step, namely to generalize Donsker's theorem to the infinite variance 
situation, is fairly obvious. It is probabely known, though we could not 
find it in the literature (compare however Bloznelis \cite{bl} and the papers 
cited therein).

\begin{proposition} Under assumption A the unconditioned processes 
$S^{n}$ converge in distribution to the L\'{e}vy-process $X$ 
fulfilling  $X_{1}$ = $\nu$ in distribution.
\end{proposition}

\emph{Proof} \quad Since $S^{n}$ and $X$ have independent, stationary 
increments, convergence of the finite-dimensional distributions 
follows immediately from assumption A. Further $|S^{n}(\tau_{n} + 
\theta_{n})- S^{n}(\tau_{n})| \stackrel{d}{=}|S^{n}(\theta_{n})| 
\rightarrow 0$ in probability for any sequence of positive numbers 
$\theta_{n} >0$, going to 0, and any sequence $\tau_{n}$ of stopping-times, 
bounded uniformly from above. Now tightness follows from a criterion 
due to Aldous \cite{al3}.  \hfill $\Box$ \bigskip

Next we discuss convergence of random walk bridges.  For the construction 
of L\'{e}vy-bridges we refer the reader to \cite{ber}, chapter VIII.3. 
Random walk bridges are treated in the next proof in quite a similar spirit.

\begin{proposition}
Under assumption A the processes $S^{n}$, conditioned on the events 
$S(n)=0$, converge in distribution to the L\'{e}vy-process $X$, 
conditioned on the event $X=0$.
\end{proposition}

\emph{Proof} \quad  Fix $t \in (0,1)$ and let $\kappa: D \rightarrow 
\mathbb{R}$ be a continuous functional, such that $\kappa(f), f \in D$ 
does only depend on the values of $f(s), s \leq t$. Then by the Markov 
property
\begin{displaymath}
	\mathrm{E}(\kappa(S^{n}) \: | S(n) =0) = \mathrm{E} (\kappa(S^{n}) 
	h_{n}(S([tn])))
\end{displaymath}
with $h_{n}(x) = \mathrm{P}(S(n)-S([tn]) = -x) / \mathrm{P}(S(n)=0)$.
By Proposition 6 
\begin{eqnarray*}
	h_{n}(x) & = & \frac{g(-x/a_{n-[tn]})}{a_{n-[tn]}} \mbox{\huge{/}}
	              \frac{g(0)}{a_{n}} (1+o(1)) \\ 
	         & = & (1-t)^{-1/\alpha}g(-x/a_{n-[tn]})/g(0)  \: (1+o(1))
\end{eqnarray*}
uniformly in $x$. Therefore
\begin{displaymath}
	\mathrm{E}(\kappa(S^{n}) \: | S(n)=0) = 
	\frac{\mathrm{E}(\kappa(S^{n})g(-(1-t)^{-1/\alpha}
	S^{n}(t)))}{(1-t)^{1/\alpha}g(0)}  (1+o(1)).
\end{displaymath}
Proposition 12 implies
\begin{displaymath}
	\mathrm{E}(\kappa(S^{n}) \: | S(n)=0) \rightarrow \frac{\mathrm{E}
	(\kappa(X)g(-(1-t)^{-1/\alpha} X(t)))}{(1-t)^{1/\alpha}g(0)}.
\end{displaymath}
This proves convergence in distribution of the process 
$(S^{n}(s))_{s \leq t }$, conditioned on the event $S(n)=0$, for any $0 <t<1$. 
Finally we have the duality relation $(S^{n}(1-s))_{s} 
\stackrel{d}{=} (S^{n}(0)-S^{n}(s-))_{s}$, thus also convergence of the 
processes $(S^{n}(s))_{s \geq 1-t}$ follows for every $0<t<1$. 
Combining these results our claim follows. \hfill $\Box$ \bigskip

The proof also shows that for any $t<1$ the distribution of the 
process $(X_{s})_{s \leq t}$, given $X_{1}=0$, is absolute continuous 
with respect to the distribution of the unconditioned 
process $(X_{s})_{s \leq t}$ (which follows as well from the 
construction in \cite{ber}).  This allows to transfer 
properties. We need the following one: The L\'{e}vy-process 
$X$ attains on $[0,1]$ its minimal value a.s. at exactly one point, namely
\begin{displaymath}
	T = \inf \{s \leq 1\: | \: X_{t} \geq \min(X_{s},X_{s-}) 
	\mbox{ for all } t \geq s\}.
\end{displaymath} 
Also $T < 1$ a.s., and $X$ is a.s. continuous at $T$. This 
follows from Propositions 2.1 to 2.4 in Millar \cite{mi} (for our process 
0 is regular for both 
$(-\infty,0)$ and $(0,\infty)$). By absolute continuity this carries 
over to the process $(X_{s})_{s\leq 1}$, conditioned on $X_{1}=0$.

There is an easy recipe using cyclic 
permutation of a path, which allows to pass over from a bridge 
$1=S(0), S(1),\ldots, S(n-1),S(n)=0$  to an 
excursion $1=S(0), S(1)>0,\ldots, S(n-1)>0,S(n)=0$, 
and which has been utilized by different people (see f.e. 
\cite{ta}). Let
\begin{displaymath}
	T^{n} = \min \{i \leq n \: | \: S(j) \geq S(i) \mbox{ for all } 
	j = i,i+1,\ldots,n \}
\end{displaymath}
be the moment, when $S$ takes its first minimal value before time $n$. Then  
$\overline{S}$, 
given by
\begin{displaymath}
	\overline{S}(i)= \left\{ \begin{array}{ll}
	 S(T^{n}+i)-S(T^{n})+1, & \mbox{ if } i \leq n-T^{n} \\
	 S(T^{n}+i-n)-S(T^{n}),  & \mbox{ if } i \geq n-T^{n},
	 \end{array} \right.
\end{displaymath}
is an excursion for any bridge $S$. In this manner we associate to each 
excursion $n+1$ different bridges, one being the given excursion itself.

This motivates to associate for any $f \in D$ an $\overline{f} \in D$, 
given by
\begin{displaymath}
	\overline{f}(t) = \left\{ \begin{array}{ll}
	f(t+ T(f)) - f(T(f)) + f(0), & \mbox{ if } t \leq 1-T(f) \\
	f(t+ T(f)-1) - f(T(f)) + f(1), & \mbox{ if } t \geq 1-T(f),
	\end{array} \right. 
\end{displaymath}
where
\begin{displaymath}
	T(f) = \min \{s \leq 1 \:| \: f(t) \geq \min(f(s),f(s-))  
	\mbox{ for all } t \geq s \}.
\end{displaymath}
Thus $\overline{f}(0)=f(0), \: \overline{f}(1)=f(1)$.

\begin{lemma}
Suppose, that $f \in D$  attains its minimum at no other point than $T(f)$. 
Then $f_{n} \rightarrow f$ implies $\overline{f}_{n} \rightarrow 
\overline{f}$ in the space $D$.
\end{lemma}

\emph{Proof} \quad $f_{n} \rightarrow f$ means that there are 
bijections $\alpha_{n}$, 
such that
\begin{displaymath}
	\sup_{t} |f_{n}(t) - f(\alpha_{n}(t))| \rightarrow 0,\; \sup_{t} 
	|\alpha_{n}(t)-t| \rightarrow 0.
\end{displaymath}
From the uniqueness of the minimum of $f$ it follows
\begin{displaymath}
	T(f_{n}) \rightarrow T(f).
\end{displaymath}
Now let $\eta_{n}>0$ be numbers going to zero and define
\begin{displaymath}
	\beta_{n}(t) = \left\{ 
	\begin{array}{ll}
		\alpha_{n}(t+T(f_{n}))-T(f), & \mbox{ if } \eta_{n} \leq t \leq 
		1-T(f_{n})  \\
		\alpha_{n}(t+T(f_{n})-1)+1-T(f), & \mbox{ if } 1-T(f_{n}) \leq t 
		\leq 1-\eta_{n}.
	\end{array} \right.
\end{displaymath}
Also put $\beta_{n}(0)=0, \: \beta_{n}(1)=1$ and continue $\beta_{n}$ 
by linear interpolation on the whole interval $[0,1]$. If 
$\eta_{n}$  goes to zero slowly enough, $\beta_{n}$ is a 
bijection of $[0,1]$, and
\begin{displaymath}
	\sup_{t} |\beta_{n}(t)-t| \rightarrow 0.
\end{displaymath}
Further
\begin{eqnarray*}
	\lefteqn{\sup_{t} |\overline{f}_{n}(t)-\overline{f}(\beta_{n}(t))|} \\
	 & \leq  & 2 \sup_{t} 
	|f_{n}(t)-f(\alpha_{n}(t))| \\
	  &  & \mbox{} +  2 \sup_{T(f_{n})-\eta_{n} \leq s,t 
	\leq T(f_{n}) + \eta_{n}} |f_{n}(s)-f(\alpha_{n}(t))| \\
	  & \leq & 4 \sup_{t} |f_{n}(t)-f(\alpha_{n}(t))| \\
	  & & \mbox{} + 2 \sup_{T(f_{n})-
	  \eta_{n} \leq s,t 
	\leq T(f_{n})+ \eta_{n}} |f(\alpha_{n}(s))-f(\alpha_{n}(t))|.
\end{eqnarray*} 
If $T(f)$ is a point of continuity of $f$, then the righthand terms 
all go to zero, and it follows the claim $\overline{f}_{n} \rightarrow 
\overline{f}$. The 
case that $f$ has a jump at $T(f)$ (which will not be considered in 
the sequel) is treated similarly; then $T(f_{n}) = T(f)$, for $n$ 
large enough. \hfill $\Box$  

\emph{Proof of Theorem 11} \quad Proposition 13  shows that the rescaled random 
walk bridges $S^{n}$ converge in distribution to a L\'{e}vy-bridge $X$. Due 
to Skorohod's theorem we may assume without loss of generality, that 
$S^{n} \rightarrow X$ a.s. in the Skorohod topology. Since $X$ 
attains its minimum a.s. at a unique point, which is a point of 
continuity, it follows from Lemma 14 $\overline{S^{n}} \rightarrow 
\overline{X}$ a.s.. By construction $Y=\overline{X}$ has the  stated 
property. \hfill $\Box$ \bigskip

This method of deriving a the normalized excursion $Y$ from the 
bridge $X$, is well-known for Brownian motion, and in the general case 
discussed in detail in Chaumont \cite{ch}.

\end{document}